\documentclass[10pt]{elsart}

\usepackage{amsmath}
\usepackage{amssymb}
\usepackage{theorem}
\usepackage{xspace}
\usepackage[a4paper]{geometry}

\def\pth#1{\left(#1\right)}
\def\acc#1{\left\{#1\right\}}
\def\cro#1{\left[#1\right]}

\def\eE{I\!\!E}
\def\eP{I\!\!P}
\def\e1{1\!\!1}

\theoremstyle{plain}

\newtheorem{remark}{Remark}[section]
\newtheorem{theorem}{Theorem}[section]
\newtheorem{corollary}{Corollary}[section]
\newtheorem{lemma}{Lemma}[section]

\newcommand{\beqn}{\begin{eqnarray*}}
\newcommand{\eeqn}{\end{eqnarray*}}

\newcommand{\N}{\mathbb{N}}
\newcommand{\R}{\mathbb{R}}

\begin{document}
\begin{frontmatter}
\title {{\bf THE S-ESTIMATOR IN CHANGE-POINT RANDOM MODEL WITH LONG MEMORY}}
\author{GABRIELA CIUPERCA }\maketitle
\address{Universit\'e de Lyon, Universit\'e Lyon 1, \\
CNRS, UMR 5208, Institut Camille Jordan, \\
Bat.  Braconnier, 43, blvd du 11 novembre 1918,\\
F - 69622 Villeurbanne Cedex, France,\\{\it email:}  Gabriela.Ciuperca@univ-lyon1.fr\\
{\it tel: }33(0)4.72.43.16.90, {\it fax: }33(0)4.72.43.16.87         
}
\begin{abstract}
The paper considers two-phase random design linear regression models. The errors and the regressors are stationary long-range dependent Gaussian. The regression parameters, the scale parameters and the change-point are estimated using a method introduced by Rousseeuw and Yohai \cite{Rousseeuw:Yohai:84}. This is called  S-estimator and it  has the property that  is more robust than the classical estimators; the outliers don't spoil the estimation results. Some asymptotic results, including the strong consistency and the convergence rate of the S-estimators, are proved.
\end{abstract}

\end{frontmatter}

\section{Introduction} 
Consider the two-phase linear regression model:
\begin{equation}
\label{eq0}
 Y_t=X_t\beta_1 \e1_{1 \leq t\leq [n\pi]}+X_t\beta_2 \e1_{[n\pi]+1 \leq t \leq n }+\varepsilon_t, \qquad t=1,...,n
\end{equation}
where $\e1_{(.)}$ is the indicator function and   $\pi \in (0,1)$, $\xi=(\beta_1,\beta_2, \pi)$, $\beta_1,\beta_2 \in \Upsilon$. The set $\Upsilon$ is  a compact  of $\R^d$, $d \geq 1$.   For this model, $Y_t$ denotes the response variable, $X_t$ is a p-vector of regressors and $\varepsilon_t$ is the error. \\
The model parameters are: regression parameters $\beta_1$ and $\beta_2$, change-point $\pi$ and error variance $\sigma^2$, with $\sigma^2 \in (0,\infty)$. Let us denote $\xi^0=(\beta^0_1,\beta^0_2,\pi^0)$ and $\sigma^2_0$  the true values of these parameters. In this paper we consider the problem of estimating of $\xi$ and $\sigma^2$, based on the observation of $(Y_t,X_t)_{1 \leq t \leq n}$. \\
Classical estimation methods studied in the statistic literature are the least squares (LS), maximum likelihood (ML) or a wider class M-estimation methods. For each of these methods one has to distinguish the cases when the errors are independent or not, and in the dependent case it is necessary to take into account the covariance structure. The same conditions can be considered for regressors $X_t$. In traditional methodology, these variables are usually assumed to be independent or with short-memory. So, if the errors are i.i.d. or with short-memory, the statistic literature related to the parametric  change-point estimation is very vast. Recent developments for the LS estimation include Feder (\cite{Feder:75a}, \cite{Feder:75b}), Bai and Perron \cite{Bai:Perron:98}, Kim and Kim \cite{Kim:Kim:08}. Bai \cite{Bai:94} considers also the least squares estimation of a shift in linear process. The process $\varepsilon_t$ is given by: $\varepsilon_t=\sum^\infty_{j=0}c_ju_{k-j} $, where $u_j$ is white noise with mean zero and variance $\sigma^2$ and the coefficients $c_j$ satisfy $\sum^\infty_{j=0}j |c_j| < \infty$. This condition excludes long-memory. For the ML estimation we refer to Bhattacharya \cite{Bhattacharya:94}, Koul and Qian \cite{Koul:Qian:02}, Ciuperca and Dapzol \cite{Ciuperca:Dapzol:08}. In the general case of the M-estimator, we can cite the papers of  Rukhin and  Vajda \cite{Rukhin:Vajda:97}, Koul et al. \cite{KQS}. Obviously, the list is not exhaustive, the subject is so large and productive that we cannot give all the papers. The convergence rate and limiting distributions of the change-point and of  the regression parameters  M-estimators are derived for the model (\ref{eq0}) by Fiteni \cite{Fiteni:02}, under restrictive and numerous assumptions. Among these conditions she considers that $(Y_t,X_t)$ is a random vector, $L_0$-NED, on a strong mixing base $\{w_t; t=...,0 ,1,... \}$, $\rho'(\varepsilon_t+\theta X_t)X_t$ is a random sequence of mean zero, $L_2$-NED of size $1/2$ on  a strong mixing base $\{w_t; t=...,0 ,1,... \}$ and $\sup_{t \leq n} \eE[ \|\rho'(\varepsilon_t+\theta X_t)X_t\|^r$ for some $r>2$. Under the same dependence assumptions, Fiteni  \cite{Fiteni:04} considers the $\tau$-estimators.\\
 On the contrary, in the case of long-memory errors or regressors, the statistical literature related to the parametric  change-point estimation is less vast.  For the simpler model:
\begin{equation}
 \label{etic}
Y_t=\mu \e1_{1 \leq t \leq k^*}+(\mu+\delta )\e1_{k^* < t \leq n}+\varepsilon_t
\end{equation}
 when the errors $\varepsilon_t$ are  long-memory Gaussian, Horvath and  Kokoszka \cite{Horvath: Kokoszka:97} considered the estimator of  $k^*$ defined by  $\hat k=\min \acc{k; |U_k|=\max_{1 \leq i <n}|U_i|}$, where, for $\gamma \in [0,1)$,  $U_i=\pth{\frac{n}{i(n-i)}}^\gamma \sum^i_{j=1}(X_j-\bar X_n)$.  The estimator converges to functionals of fractional Brownian motion. For the same model, Hidalgo and Robinson \cite{Hidalgo:Robinson:96}, Sibbertsen \cite{Sibbertsen:04} consider the LS estimator of $k$, $\mu$ and $ \delta$. A more complex model: $Y_t=\mu+(\beta+\delta \e1_{t <[\tau n]})X_t+\varepsilon_t$ is considered in the paper of Lazarov\'a \cite{Lazarova:05}, but with the supposition that $\tau$ fixed. The limiting distribution of the LS estimator of the parameters $\beta$ and $\delta$ is given.  \\
Concerning now the estimation method, it is well known that one outlier may cause a large error in a LS-estimator, ML-estimator or more generally in classical M-estimator. Nunes et al. \cite{Nunes:Kuan:Newbold:95},  Kuan and Hsu \cite{Kuan:Hsu:98} observed that, for the data that have long-memory, the LS-estimator may suggest a spurious change-point when there is none.  In that case, the parameters of the model can be estimated by using least absolute deviations (LAD) method.  If the errors are independent, Bai \cite{Bai:98} studies the LAD estimator for a multiple regime linear regression and Ciuperca \cite{Ciuperca:09} for a nonlinear change-point model. A more robust estimator was introduced by  Rousseeuw and Yohai \cite{Rousseeuw:Yohai:84}, by defining the S-estimators as the minimizers of a M-estimator of the  residual scale. The interest of the S-estimator in respect to the  LAD-estimators, is the breakdown point introduced by Hampel \cite{Hampel:71}. The breakdown point amounts to determining the smallest contaminating mass that can cause the estimator to take on value arbitrary for from the true value. Instead, concerning this method, to the author's knowledge, the past papers treat only regression models without change-point. For a linear regression model, Davies \cite{Davies:90} proves the consistency and weak convergence of S-estimator under the assumption that the errors are i.i.d. random variables.  The asymptotic behaviour of the S-estimator in a linear regression, without change-point can be also found in the papers  Zhengyan et al. \cite{Zhengyan:Degui:Jia:07}, Roeland et al. \cite{Roelant:VanAelst:Croux:08}. 
\\
In the present paper, we consider a linear regression model with a change-point in an unknown point. The regressors and the errors are assumed to be Gaussian vectors, and respectively variables, with long-memory. The regression parameters, the scale parameter and the change-point location are estimated by the S-method. The difficulty of study of the asymptotic properties of these S-estimators comes especially  from the  dependence on  change-point  in the expression of the scale parameter estimators. We first prove that the estimators are strongly convergent and afterwards their convergence rates are obtained. These rates depend of covariance structure of $X_t$ and  $\varepsilon_t$ and of Hermite rank of $\rho(\varepsilon_t/\sigma_0)-\eE[\rho(\varepsilon_t/\sigma_0)]$, where $\rho$ is the function used to construct the S-estimator. For the regression parameters and the scale parameter, we obtain the same convergence rate as in a model without change-point, let us denote it $v_n$. The S-estimator of the change-point has a faster convergence rate, more precisely $n^{-1}v_n$. This result is totally different from those obtained in the other papers where the dependence between observations is considered. Especially, let us notice that our change-point S-estimator converges more quickly towards true value than in the independence or the short-memory case or a classic  estimation method.  \\
The plan of this  paper is  as follows. In Section 2 we make some notations and assumptions afterwards we define the S-estimator for a model with change-point. In Section 3, the asymptotic behaviour of these estimator is studied. The proofs of theorems are given in  Section 4. Finally, Section 5 contains some lemmas which are useful to prove the main results.  

\section{Notation and assumptions}
Long-memory (long-range dependent) processes arise in numerous physical and social sciences. For several examples, see e.g. Baillie \cite{Baillie:96}, Cheung \cite{Cheung:93}, Lo \cite{Lo:91} among others. We  also mention  Guo and Kuol \cite{Guo: Koul:07}, where some currency exchange data sets with long-memory are considered. Another long-memory example in economy we find also in Ding et al. \cite{Ding:Granger:Engle:93} on S$\&$P daily 500 stock market returns. In that paper, they found that although the returns themselves contain little serial correlation, the absolute value of returns has significantly positive serial correlation up to 2700 lags. \\   

For the construction of the S-estimators, a function  $\rho:\R \longrightarrow [0,1]$ is needed.
 Throughout our article, we assume  that the following  classic conditions  are satisfied by  $\rho$: \\
$\bullet $ $\rho$ is symmetric, continuously  differentiable on $\R$ and $\rho(0)=0$.\\
$\bullet $  $\rho$ is increasing  in  $[0,c)$, for some $c>0$, and  constant  in $[c,\infty)$. \\
Let us denote: $\psi(z)=\rho'(z)$.\\
An example of $\rho$ satisfying these conditions was proposed by  Beaton and Tukey \cite{Beaton:Tukey:74}, for some $c>0$:
\begin{equation}
  \label{BT}
\rho(x)=\left\{ 
\begin{array}{cl}
 3(x/c)^2-3(x/c)^4+(x/c)^6, & \textrm{ if } |x| \leq c\\
1, & \textrm{ if }|x| >c
\end{array}
\right.
\end{equation}

 For model (\ref{eq0}), the following assumptions are considered: \\
{\bf (A1)} $X_t$  is a sequence of $d$-dimensional stationary long-range dependent Gaussian vectors, with  $\eE[X_t]=0$,  covariance matrix $\Gamma(t)=\eE[X_1X_{t+1}]={\cal L}(t)^TN(t) {\cal L}(t)$, where $N(t)=diag(t^{-\theta_1},...,t^{-\theta_d})$, $\theta_1,...,\theta_d \in (0,1)$ for $t \geq 1$ and $\Gamma(0)=Var(X_1)$. ${\cal L}(x)$ a $d \times d$  orthogonal matrix of slowly varying functions;\\
{\bf (A2)} $\varepsilon_t$ a sequence of stationary  long-range dependent Gaussian variables, with  $\eE[\varepsilon_t]=0$,  $\gamma(0)=Var[\varepsilon_t]=\sigma^2_0$  and the covariance $\gamma(t)=\eE[\varepsilon_1\varepsilon_{t+1}]=t^{-\alpha}L(t)$,  $\alpha \in (0,1)$ for $t \geq 1$. $L(x)$ a  positive  slowly varying function;\\
{\bf (A3)} the errors  $\varepsilon_t$ are independent of  $X_t$.\\

\noindent  The values of $\theta_1,...,\theta_d$, $\alpha$ and the functions expressions of ${\cal L}(x)$, $L(x)$ are known.\\
Recall  that a positive measurable function  $h$ is  slowly varying in Karamata's sense  if and only if, for any $\lambda>0$, ${h(\lambda x)}/{h(x)}$ converges to 1 as $x$ tends to infinity. Examples of slowly varying functions: $\log x$, $\log \log x$, $\log \log \log x$.....\\
Interested readers are referred to Beran \cite{Beran:94} or Robinson \cite{Robinson:94} for a complete reference on long-memory processes. \\
An example of process $X_t=(X_{t1},...,X_{td})$ is obtained when, for some $0< d_1 < 1/2$:
$$
X_{tj}=\sum^d_{l=1} \sum_{v \in Z} B_{jl}(t-v)\varsigma_{v,j}, \quad B_{jl}(v)=v^{-(1-d_1)}L_{jl}(v), \quad v \geq 1, \quad j,l=1,...,d
$$
where $L_{jl}$ are slowly varying functions and where $\varsigma_v=(\varsigma_{v,1},...,\varsigma_{v,d})^T$, $v \in Z$ are i.i.d. with $\varsigma_{v,j}$, $j=1,...,d$ standard Gaussian variables (see Koul and Baillie \cite{Koul:Baillie:03}). \\
For the residual function, let us consider classical notation  $r_t(\beta)=Y_t-X_t\beta$  and let $K$ the  constant  given by $K=\eE_\Phi[\rho(\varepsilon_1/\sigma_0)]$, where $\Phi$ is standard Gaussian distribution.\\
In order to construct  the S-estimator in a change-point model (\ref{eq0}), we proceed as follows:\\
- first, for $(\beta_1,\beta_2,\pi) \in \Upsilon \times \Upsilon \times (0,1)$ fixed,  scale parameter $\sigma$ is estimated by the positive solution $s_n(\xi)=s_n(\beta_1,\beta_2,\pi)$ of the equation:
\begin{equation}
 \label{eq1}
n^{-1} \sum ^{[n\pi]}_{t=1}\rho \pth{\frac{r_t(\beta_1)}{s_n(\xi)}} +n^{-1} \sum ^{n}_{t=[n\pi]+1}\rho \pth{\frac{r_t(\beta_2)}{s_n(\xi)}}=K
\end{equation}
- at the second stage, the regression parameters are estimated by the argument of the minimum of solution $s_n(\xi)$  obtained of the previous phase: 
\begin{equation}
\label{eq2}
\pth{\tilde \beta_{1n}(\pi),\tilde \beta_{2n}(\pi)}=\arg \min_{(\beta_1,\beta_2) \in \Upsilon \times \Upsilon}s_n(\beta_1,\beta_2,\pi)
\end{equation}
- in the end, the  change-point is estimated by:
\begin{equation}
 \label{eq3}
\hat \pi_n=\arg \min_{\pi \in [0,1]} s_n\pth{\tilde \beta_{1n}(\pi),\tilde \beta_{2n}(\pi),\pi}
\end{equation}
We shall make the usual identifiability assumption that the two  segments are different:
\begin{equation}
 \label{e2.2}
\beta_1 \neq \beta_2, \qquad \forall \xi \in \Upsilon \times \Upsilon \times (0,1)
\end{equation}
i.e. at least one of the  coefficients of $X_t$ has a shift. Thus the jump at $\pi$ is non-zero. This condition implies that the solution of (\ref{eq3}) is unique and it will be essential in the proof of the strong consistency. \\
 If solution $s_n(\xi) $ to (\ref{eq1}) exists then it  is  well-defined,  bounded, strictly positive, with a probability arbitrarily large (see Lemma \ref{lemma 1}).  These results are valid regardless of the covariance  structure of $X_t$, of  $\varepsilon_t$ and their distribution.   What matters is  their average is worth 0 and their variance is bounded.\\
If (\ref{eq1}) has more than one solution, $s_n(\xi)$ is defined as the supremum of all solutions. Obviously, if function $\rho$ is given by (\ref{BT}), thus equation (\ref{eq1}) has at least a solution. \\
  
In this context, we define  $\hat \sigma_n=s_n\pth{\tilde \beta_{1n}(\hat \pi_n),\tilde \beta_{2n}(\hat \pi_n),\hat \pi_n}$ as the  S-estimator of $\sigma$ and $(\hat \beta_{1n}, \hat \beta_{2n})=(\tilde \beta_{1n}(\hat \pi_n),\tilde \beta_{2n}(\hat \pi_n))$ that of $(\beta_1,\beta_2)$. We shall study the asymptotic behaviour of $\hat \sigma_n$, $(\hat \beta_{1n}, \hat \beta_{2n})$ and of $\hat \pi_n$, in the case that equation (\ref{eq1}) has at least a solution.\\

For any  $\varphi $ twice differentiable function, for $x, h \in \R$, throughout this paper we are going to use  the mean value theorem under the form:
\begin{equation}
 \label{alpha}
\varphi(x+h)=\varphi(x)+h \cro{\varphi'(x)+h \int^1_0 (1-s) \varphi'' (x+sh)ds}
\end{equation} 
For a vector $V=(v_1,\cdots, v_m)$, let us denote by $\|V\|$ its Euclidean norm and we make the convention that $|V|=(|v_1|,\cdots, |v_m|)$.\\
 In the following, we denote by $C$ a generic positive finite constant that may be different in different context, but will never depend on $n$.

\section{Asymptotic behaviour}

This section establishes asymptotic properties of the S-estimator in model (\ref{eq0}). For this purpose, first let us calculate, for solution $s_n(\xi)$ of equation  (\ref{eq1}),  the partial  derivatives  with respect to $\beta_1$ and $\beta_2$. Differentiating  (\ref{eq1}) with respect to $\beta_1$, we obtain:
$$
 \sum ^{[n\pi]}_{t=1} \frac{r_t(\beta_1)}{s_n(\xi)} \frac{\partial s_n(\xi)}{\partial \beta_1} \psi\pth{\frac{r_t(\beta_1)}{s_n(\xi)}}+  \sum^{n}_{t=[n\pi]+1}\frac{r_t(\beta_2)}{s_n(\xi)} \frac{\partial s_n(\xi)}{\partial \beta_1} \psi\pth{\frac{r_t(\beta_2)}{s_n(\xi)}} + \sum ^{[n\pi]}_{t=1} \frac{X_t}{s_n(\xi)}\psi\pth{\frac{r_t(\beta_1)}{s_n(\xi)}}=0
$$
Considering the following notation:
\begin{equation}
 \label{e5}
D_n(\xi)=n^{-1}\sum ^{[n\pi]}_{t=1} \frac{r_t(\beta_1)}{s_n(\xi)} \psi\pth{\frac{r_t(\beta_1)}{s_n(\xi)}}+n^{-1} \sum ^{n}_{t=[n\pi]+1}\frac{r_t(\beta_2)}{s_n(\xi)} \psi\pth{\frac{r_t(\beta_2)}{s_n(\xi)}}
\end{equation}
and by making similar calculation for $\partial s_n(\xi)/ \partial \beta_2$, we obtain:
\begin{equation}
\label{e6}
\left\{
\begin{array}{ll}
 \frac{\partial s_n(\xi)}{ \partial \beta_1}& =- n^{-1} D_n(\xi)^{-1}  \sum ^{[n\pi]}_{t=1}  \frac{X_t}{s_n(\xi)}\psi\pth{\frac{r_t(\beta_1)}{s_n(\xi)}}\\
\frac{\partial s_n(\xi)}{ \partial \beta_2}& = -n^{-1}  D_n(\xi)^{-1}  \sum ^{n}_{t=[n\pi]+1}
\frac{X_t}{s_n(\xi)}\psi\pth{\frac{r_t(\beta_2)}{s_n(\xi)}}
\end{array}
\right.
\end{equation}
Since $\rho$ is symmetric and increasing in $[0,c)$( and choosing suitably $c$) we have:
\begin{equation}
\label{H2}
x \psi(x) \left\{
\begin{array}{ll}
 >0, & \textrm{ if } x \in (-c,c)\setminus\{0\} \\
=0, & \textrm{ if } x=0 \textrm{ or } |x|\geq c
\end{array}
\right.
\end{equation}
By means of Lemma \ref{lemma lDn}, we prove that the random process  $D_n(\xi)^{-1}$ is bounded with a probability close to 1. In fact, the covariance structure of  $X_t$ and   of  $\varepsilon_t$, respectively, plays no role in this result. Moreover, if both random variables are no more Gaussian, Lemma \ref{lemma lDn} holds if $X_t$ and $\varepsilon_t$ are bounded with a probability close to 1.\\

In order to prove the consistency we require  that  function $\psi$ also is differentiable and strictly increasing on $(0,c)$. This condition will be used for the Taylor's expansion of $\rho$, around $(\beta^0_1,\beta^0_2)$, up to second order. \\
{\bf (H1)} $\psi(.)$ is differentiable and $\psi'(u)>0$, $\forall u \in (0,c)$.\\ 
\begin{theorem}
 \label{theorem 1}
Under assumptions  (A1)-(A3), (H1), (\ref{e2.2}), we have that    estimator $\hat \xi_n=(\hat \beta_{1n},\hat \beta_{2n}, \hat \pi_n)$ is strongly consistent: $ \hat \xi_n \overset{a.s.} {\underset{n \rightarrow \infty}{\longrightarrow}} \xi_0$. 
\end{theorem}

\begin{remark}
Statement of Theorem \ref{theorem 1} remains valid, if $X_t$ is not Gaussian,   but  it is i.i.d. and $\eE[X_tX_t^T] <\infty$. If  $\varepsilon_t$ is not Gaussian,  it has to be bounded with a probability close to 1.
\end{remark}

As a consequence of relation (\ref{e6}), the first two stages (\ref{eq1}) and (\ref{eq2}) in the construction of  the parameters  estimators, are the  solutions to the equations system:
\begin{equation}
 \label{e14}
\begin{array}{ll}
(a) & n^{-1}\sum^{[n \pi]}_{t=1} \rho \pth{\frac{r_t(\beta_1)}{\sigma}}+n^{-1}\sum^{n}_{t=[n \pi]+1} \rho \pth{\frac{r_t(\beta_2)}{\sigma}}-K=0\\
(b) &  n^{-1}\sum^{[n \pi]}_{t=1} \psi \pth{\frac{r_t(\beta_1)}{\sigma}}X_t=0\\
(c) &  n^{-1}\sum^{n}_{t=[n \pi]+1} \psi \pth{\frac{r_t(\beta_2)}{\sigma}}X_t=0
\end{array}
\end{equation}
Since the change-point intervention is essential, the convergence study of the  scale parameter estimator is realized  separately. 
According to  Theorem \ref{theorem 1}, we fix $\pi$ in a neighbourhood ${\cal V}(\pi^0)$ of $\pi^0$. In order to show the convergence of the scale parameter estimator, supplementary assumptions are needed.\\

\noindent 
{\bf (H2)} $\psi$ is twice differentiable with bounded second derivative.\\
{\bf (H3)} $\psi(x)/x$ is nonincreasing for $x>0$.\\
 Obviously, function (\ref{BT}) satisfies assumptions (H1)-(H3). 
As will be seen below,  assumption (H2) is needed to control the rest in the Taylor's expansion of $s_n(\xi)$, while (H3) is used in order to apply results  of Zhengyan et al. \cite{Zhengyan:Degui:Jia:07} on the  consistency of the scale  S-estimator in a model without change-point. Moreover, in the paper of Zhengyan et al. \cite{Zhengyan:Degui:Jia:07}, the assumption (H3) is needed to show the convergence of the regression parameter estimator, which is not the case here.  \\

\begin{theorem}
 \label{theorem 2} Under (A1)-(A3), (H1)-(H3), (\ref{e2.2}), for all  $\pi$ in a neighbourhood  ${\cal V}(\pi^0)$ of $\pi^0$, the estimator of $\sigma_0$ is strongly  consistent:
$s_n(\tilde \beta_{1n}(\pi),\tilde \beta_{2n}(\pi),\pi) \overset{a.s.} {\underset{n \rightarrow \infty}{\longrightarrow}} \sigma_0$.
\end{theorem}

\begin{corollary} 
 Under (A1)-(A3), (H1)-(H3), (\ref{e2.2}),  scale parameter S-estimator $\hat \sigma_n=s_n(\hat \beta_{1n}, \hat \beta_{2n}, \hat \pi_n)$ is strongly consistent for $\sigma_0$.
\end{corollary}

\begin{remark}
 In a model without change-point, the assumption (H2) is needed for found the convergence rate and the asymptotic distribution of the estimators but not in the convergence proof. 
\end{remark}
\begin{remark}
 The convergence result of Theorem \ref{theorem 2} holds if random vector $X_t$ is not more Gaussian but i.i.d. with $\eE[X_t]=0$ and $\eE[X_tX_t^T] < \infty$.
\end{remark}

In order  to find the convergence rate, we will use the Hermite expansion for a function of standard Gaussian variable (for details about the Hermite expansion see for example Palma \cite{Palma:07}).  Let us consider function $\chi(.):=\rho(.)-K$, where $K=\eE_\Phi[\rho(\varepsilon_1/\sigma_0)]$. Suppose  that the Hermite rank of $\chi\pth{\frac{\varepsilon_1}{\sigma_0}}$ is $q_1$. Because function $\rho$ is symmetric and $\rho(0)=0$, we have $q_1 \geq 2$.  If we denote $\nu_t=\varepsilon_t/\sigma_0$, then:
$$
\chi(\nu_t)=\sum_{q \geq {q_1}} \frac{J_q(\chi)}{q!}H_q(\nu_t)
$$
with $H_q$ the Hermite polynomial,   $J_q(\chi)=\eE[\chi(\nu_1)H_q(\nu_1)]$ and for all $ t,t'=1,\cdots n $:
\begin{equation}
\label{EHpq}
\eE[H_p(\nu_t)  H_q(\nu_{t'})]=
q! \gamma^q(t-t') \e1_{p=q}
\end{equation}
Let also $k=\min \{(\alpha q_1)/2, (\theta_i+\alpha)/2, 1 \leq i \leq d \}$.\\
In order to have the rate of convergence of the estimators in a model without change-point, following assumptions are imposed by Zhengyan et al \cite{Zhengyan:Degui:Jia:07}: $\alpha q_1 <1$ and $\max\{\alpha+\theta_j; 1 \leq j \leq d \}$.

The following theorem gives the convergence rate of  the regression parameters  and of the  scale parameter estimators. These rates are the same that in a model without change-point. \\

\begin{theorem}
 \label{theorem 3} For all $\pi \in (0,1)$, if (A1)-(A3), (H1)-(H3), (\ref{e2.2}) hold, we have
$$
\| \tilde \beta_{1n}(\pi)-\beta_{1}^0 \|=O_{\eP}\pth{(n \pi)^{-k}L_1(n\pi)}=O_{\eP}\pth{n^{-k} \tilde L_1(n)}
$$
$$
\| \tilde \beta_{2n}(\pi)-\beta_{2}^0 \|=O_{\eP}\pth{(n (1-\pi))^{-k}L_1(n(1-\pi))}=O_{\eP}\pth{n^{-k} \tilde L_1(n)}
$$
where $L_1$ and $\tilde L_1$ are slowly varying functions. For the scale parameter,  putting $\tilde s_n(\pi):=s_n(\tilde \beta_{1n}(\pi),\tilde \beta_{2n}(\pi),\pi)$, we have $|\tilde s_n(\pi) - \sigma^0|=O_{\eP}\pth{n^{-k} \tilde L_1(n)}$.
\end{theorem}

Now let us study the convergence rate of the change-point estimator: 
$$
\hat \pi_n=\arg \min_\pi s_n\pth{\tilde \beta_{1n}(\pi),\tilde \beta_{2n}(\pi),\pi}=\arg \min_\pi \cro{s_n\pth{\tilde \beta_{1n}(\pi),\tilde \beta_{2n}(\pi),\pi}-s_n(\beta_1^0,\beta^0_2,\pi)}
$$
For that we consider one of the last two equations of (\ref{e14}), for instance (c): 
$$
n^{-1} \sum^{n}_{t=[n \pi]+1} \psi \pth{\frac{r_t(\tilde \beta_{2n}(\pi))}{s_n\pth{\tilde \beta_{1n}(\pi),\tilde \beta_{2n}(\pi),\pi}}}X_t=0
$$
\begin{theorem}
\label{theorem 4} Under assumptions (A1)-(A3), (H1)-(H3), (\ref{e2.2}), we have
 $\hat \pi_n-\pi^0=O_{\eP}(n^{-1-k} \tilde L_1(n))$ with $\tilde L_1(n)$ a slowly varying function.
\end{theorem}

\noindent {\bf Example. }  If $\alpha \geq \max_{i=1,...,d}\theta_i$, then $k=(\alpha+\min_{i=1,...,d}\theta_i)/2 \leq \alpha$.\\

What is remarkable comparatively to the independence or the short-memory case is that $\hat \pi_n$  converges faster  towards $\pi_0$ when $X_t$ or $\varepsilon_t$ are long-range dependent. Consider the particular case $\alpha=\theta_1=...=\theta_d$, then $k=\alpha$.  Further  if $\alpha \in (1/2,1)$, then, for the estimators of $\beta_1$ and $\beta_2$, we have a faster convergence rate than in the independence or short-memory case. Finally, the long-memory brings about that the true values of the parameters are faster approached.\\ 
Remark also that the obtained  convergence rate completely differs from that of change-point $\tau$-estimators when  $X_t$ are   NED-dependent (Fiteni \cite{Fiteni:02}, \cite{Fiteni:04}). If $X_t$ are independent, the convergence rate is $n^{-1}$ for the change-point estimator and $n^{-1/2}$ for the parameters regression estimator, indifferently of used method: M-method (Koul et al. \cite{KQS}), ML-method (Ciuperca and Dapzol \cite{Ciuperca:Dapzol:08}), LS-estimation (Bai and Perron \cite{Bai:Perron:98}). Same convergence rate, $n^{-1}$, is obtained for change-point LS-estimator in a model with correlated errors, but not with long-memory (Bai \cite{Bai:94}).\\
It is interesting to note that the rate convergence of the change-point estimator in the mean of Gaussian variable (\ref{etic}), having long-range dependence, considered by Horvath and  Kokoszka \cite{Horvath: Kokoszka:97}, is $n^{-1}g^{-1}(1/\delta)$ with $g$ a regular varying function. Thus, the estimator of Horvath and  Kokoszka \cite{Horvath: Kokoszka:97} is slower than our estimator.\\
On the other hand, let us remark that convergence rate of the S-estimators depends of the Hermite rank of $\rho(\varepsilon_1/\sigma_0)-K$ and of the covariance structure of $X_t$ and $\varepsilon_t$. 

\section{Proofs of Theorems}
{\bf Proof of Theorem  \ref{theorem 1}.} Let us consider the function $e(\xi)= \eE [s_n(\eta,\pi)-s_n(\eta^0,\pi^0)]$, with supposition, without loss the generality,  that $\pi \leq \pi^0$.  Using the same arguments as for  (\ref{e11}), we obtain that:  $
\eE[|s_n(\eta,\pi)-s_n(\eta,\pi^0)|] \leq C \| \beta_1- \beta_2\| \cdot |\pi-\pi^0| <\infty 
$ 
and similarly to  (\ref{e8}): $
\eE[|s_n(\eta,\pi^0)-s_n(\eta^0,\pi^0)|] \leq C \| \eta-\eta^0 \|$. 
Thus, function $e(\xi)$ is well-defined. By  Lemma \ref{lemma 2},  function $e(\xi)$ is  continuous and furthermore $e(\xi^0)=0$. 
For using an argument like the one in Huber \cite{Huber:67}, we will to prove that:  $\eE[s_n(\eta,\pi) - s_n(\eta^0,\pi^0)]>0$, for every $ \xi \neq \xi^0$. \
Since $s_n(\xi)$ and $s_n(\xi^0)$ are both solutions of equation (\ref{eq1}), we have $0 = (S_{1,n}^{(0)}+S_{1,n}^{(1)})+(S_{2,n}^{(0)}+S_{2,n}^{(1)})+(S_{3,n}^{(0)}+S_{3,n}^{(1)})$, with:
$$
S_{1,n}^{(0)} \equiv n^{-1} \sum^{[ n\pi]}_{t=1} \cro{ \rho \pth{\frac{r_t(\beta_1^0)}{s_n(\xi)}}- \rho \pth{\frac{r_t(\beta_1^0)}{s_n(\xi^0)}}}, \; S_{1,n}^{(1)} \equiv n^{-1} \sum^{[ n\pi]}_{t=1} \cro{\rho \pth{\frac{r_t(\beta_1)}{s_n(\xi)}}-\rho \pth{\frac{r_t(\beta_1^0)}{s_n(\xi)}}}
$$
$$
S_{2,n}^{(0)} \equiv n^{-1} \sum^{[n \pi^0]}_{t=[n \pi ]+1} \cro{\rho \pth{\frac{r_t(\beta_1^0)}{s_n(\xi)}}- \rho \pth{\frac{r_t(\beta_1^0)}{s_n(\xi^0)}}}, \; S_{2,n}^{(1)} \equiv n^{-1} \sum^{[n \pi^0]}_{t=[n \pi ]+1} \cro{\rho \pth{\frac{r_t(\beta_2)}{s_n(\xi)}}-\rho \pth{\frac{r_t(\beta_1^0)}{s_n(\xi)}} }
$$
$$
S_{3,n}^{(0)} \equiv n^{-1} \sum^{n}_{t=[n \pi^0]+1} \cro{ \rho \pth{\frac{r_t(\beta_2^0)}{s_n(\xi)}}- \rho \pth{\frac{r_t(\beta_2^0)}{s_n(\xi^0)}} }, \; S_{3,n}^{(1)} \equiv n^{-1} \sum^{n}_{t=[n \pi^0]+1} \cro{\rho \pth{\frac{r_t(\beta_2)}{s_n(\xi)}}-\rho \pth{\frac{r_t(\beta_2^0)}{s_n(\xi)}} }
$$
Then, by the mean value theorem (TVM), $S_{1,n}^{(0)}+S_{2,n}^{(0)}+S_{3,n}^{(0)}$ can be written as:
$$
n^{-1}\pth{\frac{1}{s_n(\xi)}-\frac{1}{s_n(\xi^0)}} \left[\sum^{[ n\pi]}_{t=1} r_t(\beta^0_1) \psi \pth{\frac{r_t(\beta^0_1)}{u^{(1)}_n(\eta^0,\pi,\pi^0)}}+\sum^{[n \pi^0]}_{t=[n \pi ]+1}r_t(\beta^0_1) \psi \pth{\frac{r_t(\beta^0_1)}{u^{(2)}_n(\eta^0,\pi,\pi^0)}} \right.$$
$$
 +\left. \sum^{n}_{t=[n \pi^0]+1} r_t(\beta^0_2) \psi \pth{\frac{r_t(\beta^0_2)}{u^{(3)}_n(\eta^0,\pi,\pi^0)}} \right]
$$
with $u^{(1)}_n,u^{(2)}_n,u^{(3)}_n$ defined in the same way as in the proof of the Lemma \ref{lemma 2}.
Moreover, using property (\ref{H2}), we have the following:  $S_{1,n}^{(0)}+S_{2,n}^{(0)}+S_{3,n}^{(0)}=[s_n(\xi^0)-s_n(\xi)] V_n$, where $V_n$ is  a positive random variable with probability close to 1. \\
Moreover, using Taylor's expansion, the expressions of   $S_{1,n}^{(1)}$, $S_{2,n}^{(1)}$  and $S_{3,n}^{(1)}$ can be written as:
$$
S_{1,n}^{(1)}=n^{-1}\sum^{[ n\pi]}_{t=1} X_t(\beta^0_1-\beta_1)\cro{\psi \pth{\frac{r_t(\beta^0_1)}{s_n(\xi)}}+\frac 12  \psi'\pth{\frac{\varepsilon_t+\delta_1X_t(\beta^0_1-\beta_1)}{s_n(\xi)}}(\beta^0_1-\beta_1)^TX^T_t}
$$ 
$$
S_{2,n}^{(1)}=n^{-1} \sum^{[n \pi^0]}_{t=[n \pi ]+1}X_t(\beta^0_1-\beta_2)\cro{\psi \pth{\frac{r_t(\beta^0_1)}{s_n(\xi)}}+\frac 12  \psi'\pth{\frac{\varepsilon_t+\delta_2 X_t(\beta^0_1-\beta_2)}{s_n(\xi)}}(\beta^0_1-\beta_2)^TX^T_t}
$$
$$
S_{3,n}^{(1)}=n^{-1} \sum^{n}_{t=[n \pi^0]+1}  X_t(\beta^0_2-\beta_2) \cro{\psi \pth{\frac{r_t(\beta^0_2)}{s_n(\xi)}}+\frac 12  \psi'\pth{\frac{\varepsilon_t+\delta_3 X_t(\beta^0_2-\beta_2)}{s_n(\xi)}}(\beta^0_2-\beta_2)^TX^T_t}
$$
with $\delta_1,\delta_2,\delta_3 \in (0,1)$. By the ergodic theorem, we  obtain:
\begin{equation}
 \label{e12}
\begin{array}{l}
n^{-1}\sum^{[ n\pi]}_{t=1} X_t(\beta^0_1-\beta_1)\psi \pth{\frac{r_t(\beta^0_1)}{s_n(\xi)}}=o_{\eP}(1), n^{-1} \sum^{[n \pi^0]}_{t=[n \pi ]+1}X_t(\beta^0_1-\beta_2)\psi \pth{\frac{r_t(\beta^0_1)}{s_n(\xi)}}=o_{\eP}(1),\\
n^{-1} \sum^{n}_{t=[n \pi^0]+1} X_t(\beta^0_2-\beta_2)\psi \pth{\frac{r_t(\beta^0_2)}{s_n(\xi)}}=o_{\eP}(1)
\end{array}
\end{equation}
Relation (\ref{e12}) and assumption (H1) imply: for  any $\xi \neq \xi^0$, for all $\epsilon >0$,  there exits $ a>0$, such that
\begin{equation}
 \label{e13}
\eP[S_{1,n}^{(1)}+S_{2,n}^{(1)}+S_{3,n}^{(1)}>a] >1-\epsilon
\end{equation}
Assumption (\ref{e2.2}), the above relation and  $S_{1,n}^{(1)}+S_{2,n}^{(1)}+S_{3,n}^{(1)}=-(S_{1,n}^{(0)}+S_{2,n}^{(0)}+S_{3,n}^{(0)})=[s_n(\xi)-s_n(\xi^0)] V_n$, with  $V_n>0$, imply the  conclusion $\eE[s_n(\eta,\pi) - s_n(\eta^0,\pi^0)]>0$, for all $\xi \neq \xi^0$. Using this, the compactness of the parameter space, $\hat \xi_n=\arg \min_{\xi \in \Upsilon \times \Upsilon \times [0,1]} s_n(\xi)$ and an argument like one in Huber \cite{Huber:67}, the strongly convergence of $\hat \xi_n$ results.
\hspace*{\fill}$\blacksquare$ \\

{\bf Proof of Theorem  \ref{theorem 2}.}
We first  prove that, if  we consider in (\ref{eq0}) the true value for $\eta$ and $\pi$, then the scale parameter estimator is strongly consistent:
\begin{equation}
 \label{e15}
s_n(\eta^0,\pi^0) \overset{a.s.} {\underset{n \rightarrow \infty}{\longrightarrow}} \sigma_0
\end{equation}
Let us observe that in fact $s_n(\eta^0,\pi^0)$ is  the  solution of a problem without breaking:
$$
K=n^{-1}\sum^{[n \pi^0]}_{t=1} \rho \pth{\frac{\varepsilon_t}{s_n(\xi^0)}}+n^{-1}\sum^{n}_{t=[n \pi^0]+1}\rho \pth{\frac{\varepsilon_t}{s_n(\xi^0)}} =n^{-1}\sum^n_{t=1} \rho \pth{\frac{\varepsilon_t}{s_n(\xi^0)}}
$$
and then, relation (\ref{e15}) is obtained  by Theorem 3.1 of Zhengyan et al.  \cite{Zhengyan:Degui:Jia:07}. Now,  as a consequence  of Theorem \ref{theorem 1}, we may  consider only the case  $(\eta,\pi)$ in a neighbourhood  ${\cal V}(\eta^0,\pi^0)$ of $(\eta^0,\pi^0)$. Consider the decomposition:
\begin{equation}
\label{dd}
 s_n(\eta,\pi)-s_n(\eta^0,\pi^0)=[s_n(\eta,\pi)-s_n(\eta^0,\pi)]+[s_n(\eta^0,\pi)-s_n(\eta^0,\pi^0)]:\equiv S_1(n)+S_2(n)
\end{equation}
Since  $S_1(n)$, depends only on the regression parameters, by Theorem  \ref{theorem 1}, taking  into account  relations (\ref{e6}) and (\ref{e7}), we readily obtain:
\begin{equation}
\label{point}
\sup_{\eta \in {\cal V}(\eta^0)} | s_n(\eta,\pi)-s_n(\eta^0,\pi) | \overset{a.s.} {\underset{n \rightarrow \infty}{\longrightarrow}} 0
\end{equation}
For $S_2(n)$, an argument like the one used for  (\ref{e10}) yield that  $s_n(\eta^0,\pi)-s_n(\eta^0,\pi^0)$  behaves as: 
$$
n^{-1} \sum^{[n \pi^0]}_{t=[n \pi]+1} X_t \psi\pth{\frac{\tilde r_t(\beta^0_1,\beta^0_2)}{s_n(\eta^0,\pi)}}
$$
where $\tilde r_t(\beta^0_1,\beta^0_2)= r_t(\beta^0_1)+m_t[r_t(\beta^0_2)-r_t(\beta^0_1)]$, with $0<m_t <1$. 
Let us remark that $\tilde r_t(\beta^0_1,\beta^0_2)=\varepsilon_t$. We write Taylor's expansion of $\psi\pth{{\varepsilon_t}/{s_n(\eta^0,\pi)}}$ around  $\psi\pth{{\varepsilon_t}/{\sigma_0}}$ up to second order:
$$
n^{-1} \sum^{[n \pi^0]}_{t=1} X_t \psi\pth{\frac{\varepsilon_t}{s_n(\eta^0,\pi)}}=n^{-1} \sum^{[n \pi^0]}_{t=1} X_t \psi\pth{\frac{\varepsilon_t}{\sigma_0}}- n^{-1}\frac{s_n(\eta^0,\pi)- \sigma_0}{\sigma_0 s_n(\eta^0,\pi)}\sum^{[n \pi^0]}_{t=1} X_t \varepsilon_t \psi'\pth{\frac{\varepsilon_t}{\sigma_0}}
$$
$$
+n^{-1}\frac{(s_n(\eta^0,\pi)- \sigma_0)^2}{2 \sigma_0 s_n(\eta^0,\pi)}\sum^{[n \pi^0]}_{t=1} \psi''(\varsigma_t)\varepsilon^2_tX_t
$$
with  $\varsigma_t=\varepsilon_t[s_n(\eta^0,\pi)+\upsilon_t(\sigma_0-s_n(\eta^0,\pi))]/(\sigma_0s_n(\eta^0,\pi))$, $\upsilon_t \in (0,1)$.  Since  $\psi''$ is bounded, we have  $n^{-1}\sum^{[n \pi^0]}_{t=1} \psi''(\varsigma_t)\varepsilon^2_tX_t < \infty $ with  probability 1. Moreover:
$$
\pi^0 \cro{\frac{1}{[n \pi^0]}\sum^{[n \pi^0]}_{t=1} X_t \varepsilon_t \psi'\pth{\frac{\varepsilon_t}{\sigma_0}} }\overset{a.s.} {\underset{n \rightarrow \infty}{\longrightarrow}}  \pi^0 \eE \cro{X_t  \varepsilon_t \psi'\pth{\frac{\varepsilon_t}{\sigma_0}}}=0
$$ 
Hence: $
n^{-1} \sum^{[n \pi^0]}_{t=1} X_t \cro{ \psi\pth{{\varepsilon_t}/{s_n(\eta^0,\pi)}} -\psi\pth{{\varepsilon_t}/{\sigma_0}} } =o_{\eP}(s_n(\eta^0,\pi)- \sigma_0)$. This relation and $
n^{-1} \sum^{[n \pi^0]}_{t=1} X_t \psi\pth{{\varepsilon_t}/{\sigma_0}} \overset{a.s.} {\underset{n \rightarrow \infty}{\longrightarrow}} 0$ yield that $S_2(n)=s_n(\eta^0,\pi)-s_n(\eta^0,\pi^0)=o_{\eP}(1)+o_{\eP}(s_n(\eta^0,\pi)-\sigma_0)=o_{\eP}(1)+o_{\eP}(S_2(n))$, for the last relation we have used (\ref{e15}). Then  $\sup_{\pi \in {\cal V}(\pi^0)}|S_2(n)|\overset{a.s.} {\underset{n \rightarrow \infty}{\longrightarrow}} 0 $. This fact, with relation  (\ref{point}), together with decomposition (\ref{dd}) and relation (\ref{e15}), yield the Theorem.
\hspace*{\fill}$\blacksquare$ \\

{\bf Proof of Theorem  \ref{theorem 3}.}
For  $\pi \in (0,1)$ fixed, the convergence rate of the  regression parameters estimator $\tilde \beta_{1n}(\pi)$ and $\tilde \beta_{2n}(\pi)$ is obtained by the application of Zhengyan et al. \cite{Zhengyan:Degui:Jia:07} results on every segment. On the other hand, the study of the convergence rate of $\hat s_n$ is more difficult because it interferes in both segments. For notational simplicity, in the rest of this proof,  we denote $\tilde \beta_{1n}=\tilde \beta_{1n}(\pi)$,  $\tilde \beta_{2n}=\tilde \beta_{2n}(\pi)$ and $\tilde s_n=\tilde s_n(\pi)$. The study will be made in two stages. First, we are going to write equation  (\ref{e14})(a) in  another form, putting in evidence $\sigma_0$ by a limited development. Afterwards, in the second stage, the obtained form is studied by taking into account the convergence rate of the regression parameters estimators and what $X_t$, $\varepsilon_t$ are long-memory Gaussian.\\
{\it Stage 1}.  
 Equation (\ref{e14})(a) can be expressed as:
\begin{equation}
 \label{e17}
n^{-1}\sum^{[n \pi]}_{t=1} \chi \pth{\frac{r_t(\tilde \beta_{1n})}{\tilde s_n}}+n^{-1}\sum_{t=[n \pi]+1}^{n} \chi \pth{\frac{r_t(\tilde \beta_{2n})}{\tilde s_n}}=0
\end{equation}
We apply  (\ref{alpha}) to  function $\chi$ and for:\\
$t=1,\cdots, [n \pi]$, $x_t={r_t(\tilde \beta_{1n})}/{\tilde s_n}$, $h_t=\pth{\sigma_0^{-1}-\tilde s_n^{-1}} r_t(\tilde \beta_{1n})$\\
$t=[n \pi]+1,\cdots, n$, $x_t={r_t(\tilde \beta_{2n})}/{\tilde s_n}$, $h_t=\pth{\sigma_0^{-1}-\tilde s_n^{-1}} r_t(\tilde \beta_{2n})$\\
Hence, for the part $t=1,\cdots, [n \pi]$,  we have: 
\begin{equation}
 \label{e18}
\begin{array}{ll}
n^{-1}\overset{[n \pi]}{\underset{t=1}{\sum}} \chi (x_t+h_t)= & n^{-1}\overset{[n \pi]}{\underset{t=1}{\sum}}\chi (x_t)+\frac{\tilde s_n - \sigma_0}{\sigma_0 \tilde s_n} \left[ n^{-1}\overset{[n \pi]}{\underset{t=1}{\sum}}  r_t(\tilde \beta_{1n}) \psi(x_t) \right.\\
& + \left.n^{-1}  \frac{\tilde s_n- \sigma_0}{\sigma_0 \tilde s_n} \overset{[n \pi]}{\underset{t=1}{\sum}}  r^2 _t(\tilde \beta_{1n}) \int^1_0 (1-s) \psi' \pth{\frac{r_t(\tilde \beta_{1n})}{\sigma_0}+s \cdot  r_t(\tilde \beta_{1n}) \pth{\frac{1}{\sigma_0}-\frac{1}{\tilde s_n}}} ds \right]
 \end{array}
\end{equation}
Thus, in order to study the first sum of (\ref{e17}), we shall analyse the terms of the right-hand side of (\ref{e18}). \\
We first consider the last term of  the right-hand side of (\ref{e18}). Elementary algebra yields that:
$$
n^{-1} \frac{\tilde s_n- \sigma_0}{\sigma_0 \tilde s_n} \sum^{[n \pi]}_{t=1} r^2_t(\tilde \beta_{1n})=\frac{\tilde s_n- \sigma_0}{\sigma_0 \tilde s_n} \left[n^{-1} \sum^{[n \pi]}_{t=1} \varepsilon^2_t +n^{-1}  \sum^{[n \pi]}_{t=1}(\tilde \beta_{1n} - \beta^0_1)^T X^T X_t (\tilde \beta_{1n} - \beta^0_1)\right.
$$
$$
+\left.  2n^{-1}  \sum^{[n \pi]}_{t=1} \varepsilon_t X_t (\tilde \beta_{1n} - \beta^0_1) \right]
$$
By the ergodic theorem $n^{-1}\sum^{[n \pi]}_{t=1}\varepsilon^2_t=O_{\eP}(1)$, $n^{-1}\sum^{[n \pi]}_{t=1}X^T X_t=O_{\eP}(1)$ and since  $\varepsilon_t$ and $X_t$ are independent, we have $n^{-1}\sum^{[n \pi]}_{t=1}\varepsilon_t X_t=o_{\eP}(1)$.   Thus, since $\psi'$ is bounded, $\tilde s_n-\sigma_0=o_{\eP}(1)$, $\tilde s_n >0$ with probability 1, the last  term of the right-hand side of  (\ref{e18}) is $o_{\eP}(1)$. \\
We now consider the second term of  the right-hand side of (\ref{e18}). For the sum, we have:
$$
n^{-1} \left\| \sum^{[n \pi]}_{t=1} r_t(\tilde \beta_{1n}) \psi \pth{\frac{r_t(\tilde \beta_{1n})}{\tilde s_n}} - \sum^{[n \pi]}_{t=1} r_t(\beta^0_1) \psi \pth{\frac{r_t(\beta^0_1)}{\sigma_0}} \right\|
$$
$$
\leq n^{-1} \left\| \sum^{[n \pi]}_{t=1} \cro{r_t(\tilde \beta_{1n}) - r_t(\beta^0_1)}\psi \pth{\frac{r_t(\tilde \beta_{1n})}{\tilde s_n}} \right\|+n^{-1} \left\| \sum^{[n \pi]}_{t=1} r_t(\beta^0_1) \cro{\psi \pth{\frac{r_t(\tilde \beta_{1n})}{\tilde s_n}} - \psi \pth{\frac{r_t(\beta^0_1)}{\sigma_0}}}  \right\|
$$
and since  $\psi$ is bounded:
$$
\leq C\|\tilde \beta_{1n} - \beta^0_1\| n^{-1}  \sum^{[n \pi]}_{t=1} \left\|X_t\right\| +C n^{-1}  \sum^{[n \pi]}_{t=1}\left\| \varepsilon_t\right\|
$$
The above inequality,  with the ergodic theorem, $\eE[\|X_t\|] < \infty$, $\eE[\|\varepsilon_t\|] < \infty$ 
and $\tilde \beta_{1n} - \beta^0_1=o_{\eP}(1)$, imply that 
$$
n^{-1}\sum^{[n \pi]}_{t=1} r_t(\tilde \beta_{1n}) \psi \pth{\frac{r_t(\tilde \beta_{1n})}{\tilde s_n}} - n^{-1}\sum^{[n \pi]}_{t=1} r_t(\beta^0_1) \psi \pth{\frac{r_t(\beta^0_1)}{\sigma_0}}=o_{\eP}(1)
$$
Thus, the second term of  the right-hand side of (\ref{e18}) can be expressed:
$$
n^{-1} \sum^{[n \pi]}_{t=1} h_t \psi(x_t)=\frac{\tilde s_n- \sigma_0}{\sigma_0 \tilde s_n} \pth{n^{-1} \sum^{[n \pi]}_{t=1}  \varepsilon_t \psi \pth{\frac{\varepsilon_t}{\sigma_0}}+o_{\eP}(1)}
$$
Then, relation (\ref{e18}) becomes:
\begin{equation}
 \label{e22}
n^{-1} \sum^{[n \pi]}_{t=1} \chi (x_t+h_t)=  n^{-1} \sum^{[n \pi]}_{t=1} \chi (x_t)+\frac{\tilde s_n- \sigma_0}{\sigma_0 \tilde s_n} \pth{n^{-1} \sum^{[n \pi]}_{t=1}  \varepsilon_t \psi \pth{\frac{\varepsilon_t}{\sigma_0}}+o_{\eP}(1)}
\end{equation}
A similar relation holds for the part  $t=[n \pi]+1, \cdots, n$:
\begin{equation}
 \label{e23}
n^{-1} \sum_{t=[n \pi]+1}^{n} \chi (x_t+h_t)=  n^{-1} \sum_{t=[n \pi]+1}^{n}  \chi (x_t)+\frac{\tilde s_n- \sigma_0}{\sigma_0 \tilde s_n} \pth{n^{-1} \sum_{t=[n \pi]+1}^{n}   \varepsilon_t \psi \pth{\frac{\varepsilon_t}{\sigma_0}}+o_{\eP}(1)}
\end{equation}
Adding  (\ref{e22}) and (\ref{e23}), taking into account the relation (\ref{e17}), we obtain:
$$
0=n^{-1} \sum^{[n \pi]}_{t=1} \chi \pth{\frac{r_t(\tilde \beta_{1n})}{\sigma_0}}+n^{-1} \sum_{t=[n \pi]+1}^{n} \chi \pth{\frac{r_t(\tilde \beta_{2n})}{\sigma_0}}+\frac{\tilde s_n- \sigma_0}{\sigma_0 \tilde s_n} \pth{n^{-1} \sum_{t=1}^{n}   \varepsilon_t \psi \pth{\frac{\varepsilon_t}{\sigma_0}}+o_{\eP}(1)}
$$
By ergodic theorem: $
n^{-1} \sum_{t=1}^{n}   \varepsilon_t \psi \pth{{\varepsilon_t}/{\sigma_0}} \overset{\eP} {\underset{n \rightarrow \infty}{\longrightarrow}} \eE \cro{\varepsilon_1 \psi \pth{{\varepsilon_1}/{\sigma_0}}}$. \\
{\it Stage 2}.  Then, the convergence rate of $\tilde s_n$ will be obtained by studying:
\begin{equation}
 \label{e24}
n^{-1} \sum^{[n \pi]}_{t=1} \chi \pth{\frac{r_t(\tilde \beta_{1n})}{\sigma_0}}+n^{-1} \sum_{t=[n \pi]+1}^{n} \chi \pth{\frac{r_t(\tilde \beta_{2n})}{\sigma_0}}= \frac{\sigma_0-\tilde s_n}{\sigma_0 \tilde s_n}\cro{\eE \cro{\varepsilon_1 \psi \pth{\frac{\varepsilon_1}{\sigma_0}}}+o_{\eP}(1)}
\end{equation}
For $t=1, \cdots, [n \pi]$,  making the Taylor's expansion of  $\chi$ up to second order, we obtain  that $n^{-1} \sum^{[n \pi]}_{t=1} \chi \pth{\sigma_0^{-1} r_t(\tilde \beta_{1n})}$ can be written as:
\begin{equation}
\label{chi}
n^{-1}\left\{\sum^{[n \pi]}_{t=1}\chi\pth{\frac{\varepsilon_t}{\sigma_0}}-\frac{1}{\sigma_0}\sum^{[n \pi]}_{t=1}\chi'\pth{\frac{\varepsilon_t}{\sigma_0}}X_t(\tilde \beta_{1n}- \beta_1^0)-\frac{1}{2 \sigma^2_0}\sum^{[n \pi]}_{t=1} \chi''\pth{\frac{\varepsilon_t-\delta_tX_t(\tilde \beta_{1n}- \beta_1^0)}{\sigma_0}}[X_t(\tilde \beta_{1n}- \beta_1^0)]^2 \right\}
\end{equation}

Let us analyse the three terms of the previous equation separately. \\
$\bullet$ For the first term, let us $\nu_t=\varepsilon_t/\sigma_0 \sim {\cal N}(0,1)$ denote. We use the Hermite  expansion for $\sum^{[n \pi]}_{t=1} \chi(\nu_t)$. Because the Hermite rank of $\chi(\nu_t)$ is $q_1$, $q_1 \geq 2$, by  (\ref{EHpq}) below:
\begin{equation}
\label{I1I2}
\sum^{[n \pi]}_{t=1} \chi(\nu_t)=\frac{J_{q_1}(\chi)}{q_1!} \sum^{[n \pi]}_{t=1} H_{q_1}(\nu_t)+\sum^{[n \pi]}_{t=1}  \sum_{q \geq q_1+1} \frac{J_q(\chi)}{q!}H_q(\nu_t):\equiv  T_{1,n}+T_{2,n} 
\end{equation}
\noindent For $T_{1,n}$ we have:
$$
\eE[T_{1,n}^2]=\frac{J^2_{q_1}(\chi)}{(q_1!)^2}  \sum^{[n \pi]}_{t=1} \sum^{[n \pi]}_{j=1} (q_1!)\gamma^{q1}(|t-j|)=(q_1!)\frac{J^2_1(\rho)}{(q_1!)^2} \cro{[n \pi] \gamma^{q_1}(0)+2 \sum^{[n \pi]-1}_{t=1} ([n \pi]-t) \gamma^{q_1}(t)}
$$
$$
=(q_1!) \frac{J^2_1(\rho)}{(q_1!)^2} \cro{O(n) +2 \sum^{[n \pi]-1}_{t=1} ([n \pi]-t) t^{-\alpha q_1} L^{q_1}(t)}= \frac{J^2_1(\rho)}{q_1!} \cro{O(n) +O(n^{2-\alpha q_1}) L^{q_1}([n\pi]) }$$
$$
=O(n^{2-\alpha q_1}) L^{q_1}([n\pi])
$$
For  $T_{2,n}$ we have:
$$
\eE[T_{2,n}^2]=\sum_{q \geq q_1+1} \frac{J_q^2(\rho)}{q!} \cro{\sum^{[n \pi]}_{t=1} \sum^{[n \pi]}_{j=1}  \gamma^q(|t-j|)} \leq \sum_{q \geq q_1+1}  \frac{J_q^2(\rho)}{q!}\sum^{[n \pi]}_{t=1} \sum^{[n \pi]}_{j=1}  \gamma^{q_1+1}(|t-j|) 
$$
$$
=O(n)+2\sum_{q \geq q_1+1}  \frac{J_q^2(\rho)}{q!}   \sum^{[n \pi]-1}_{t=1} ([n \pi]-t) \gamma^{q_1+1}(t) \leq O(n)+ 2\sum_{q \geq q_1+1}  \frac{J_q^2(\rho)}{q!}   \sum^{[n \pi]-1}_{t=1} ([n \pi]-t) 
t^{-(q_1+1) \alpha } L^{q_1+1}(t)$$
$$ =O(n^{2-(q_1+1) \alpha} L^{q_1+1}([n \pi]))
$$
Hence $\eE[T_{2,n}^2]=o(\eE[T_{1,n}^2])$. Then, for equation (\ref{I1I2}), we straightforwardly have: 
\begin{equation}
\label{t1}
\sum^{[n \pi]}_{t=1} \chi(\nu_t)=O_{\eP} \pth{{\eE[T_{1,n}^2]}}^{1/2}=O_{\eP}(n^{1-\alpha q_1 / 2}) L^{q_1/2}([n\pi])
\end{equation}
$\bullet$ For the second term of  (\ref{chi}), since $\nu_t$ and $X_t$ are independent, by ergodic theorem,  we have:
\begin{equation}
 \label{t2}
n^{-1} \sum^{[n \pi]}_{t=1}\chi'\pth{\nu_t}X_t(\tilde \beta_{1n}- \beta_1^0)=o_{\eP}(\|\tilde \beta_{1n}- \beta_1^0 \|)
\end{equation}
$\bullet$ For the third term of  (\ref{chi}), since $\psi'$ is bounded and $n^{-1} \sum^{[n \pi]}_{t=1} X_tX_t^T=O_{\eP}(1)$, we have:
\begin{equation}
 \label{t3}
n^{-1} \sum^{[n \pi]}_{t=1} \chi''\pth{\nu_t}[X_t(\tilde \beta_{1n}- \beta_1^0)]^2=O_{\eP}(\|\tilde \beta_{1n}- \beta_1^0 \|^2)=o_{\eP}(\|\tilde \beta_{1n}- \beta_1^0 \|)
\end{equation}
Then, by taking  (\ref{t1}), (\ref{t2}), (\ref{t3}) into account, the behaviour of  (\ref{chi}) is given by (\ref{t1}) and it is  $O_{\eP}(n^{-\alpha q_1/ 2}) L^{q_1/2}([n\pi])+o_{\eP}(\| \tilde \beta_{1n}-\beta^0_1\|)$. Similar one reasoning is made for the part   $t=[n \pi]+1,\cdots,n$ and we obtain that: $n^{-1}\sum^n_{t=[n \pi]+1}\chi \pth{\sigma_0^{-1} r_t(\tilde \beta_{2n})}=O_{\eP}(n^{-\alpha q_1/ 2}) L^{q_1/2}(n (1-[\pi]))+o_{\eP}(\| \tilde \beta_{2n}-\beta^0_2\|)$. Then, for  relation (\ref{e24}), we have:
$$
\frac{\sigma_0-\tilde s_n}{\sigma_0 \tilde s_n}\cro{\eE \cro{\varepsilon_1 \psi \pth{\frac{\varepsilon_1}{\sigma_0}}}+o_{\eP}(1)}=O_{\eP}(n^{-\alpha q_1/ 2}) L^{q_1/2}(n)+o_{\eP}(\| \tilde \beta_{1n}-\beta^0_1\|+\| \tilde \beta_{2n}-\beta^0_2\|)
$$
and the convergence rate of $\tilde s_n$ follows.
\hspace*{\fill}$\blacksquare$ \\

{\bf Proof of Theorem  \ref{theorem 4}.} As a consequence of  Theorem \ref{theorem 1}, we consider $\pi$ in a neighbourhood of $\pi^0$. We suppose, without loss of generality, that $\pi <\pi^0$. 
Considering relation  (\ref{e14})(c), we have:
\begin{equation}
 \label{e25}
n^{-1} \sum^{[n \pi^0]}_{t=[n \pi]+1} \psi \pth{\frac{r_t(\tilde \beta_{2n}(\pi))}{s_n\pth{\tilde \beta_{1n}(\pi),\tilde \beta_{2n}(\pi),\pi}}}X_t=-n^{-1}\sum^{n}_{t=[n \pi^0]+1} \psi \pth{\frac{r_t(\tilde \beta_{2n}(\pi))}{s_n\pth{\tilde \beta_{1n}(\pi),\tilde \beta_{2n}(\pi),\pi}}}X_t
\end{equation}
Since $\| \tilde \beta_{2n}(\pi)- \beta^0_2\|=O_{\eP}(n^{-k}\tilde L_1(n))$, an argument like the one used for relation (\ref{t2}) yield that the right-hand side of (\ref{e25}) is $O_{\eP}(n^{-k}\tilde L_1(n))$.\\
We apply  (\ref{alpha}) to function $\psi$, for: $x_t=\frac{\varepsilon_t}{s_n\pth{\tilde \beta_{1n}(\pi),\tilde \beta_{2n}(\pi),\pi} }$, $h_t=-\frac{X_t(\tilde \beta_{2n}(\pi) - \beta^0_1)}{s_n\pth{\tilde \beta_{1n}(\pi),\tilde \beta_{2n}(\pi),\pi}}$. 
For the left-hand side of  (\ref{e25}),  since  $s_n\pth{\tilde \beta_{1n}(\pi),\tilde \beta_{2n}(\pi),\pi} \rightarrow \sigma_0$ a.s. for $n \rightarrow \infty $, and $\beta^0_1 \neq \beta^0_2$,  we obtain :
\begin{equation}
 \label{e26}
n^{-1} \sum^{[n \pi^0]}_{t=[n \pi]+1} \psi \pth{\frac{r_t(\tilde \beta_{2n}(\pi))}{s_n\pth{\tilde \beta_{1n}(\pi),\tilde \beta_{2n}(\pi),\pi}}}X_t=n^{-1}\sum^{[n \pi^0]}_{t=[n \pi]+1} \psi \pth{\frac{\varepsilon_t}{\sigma_0}}X_t+O_{\eP}(n(\pi^0-\pi))
\end{equation}
But, making  Hermite expansion of  $\psi(\nu_t)$, we get:
$$
\sum^{[n \pi^0]}_{t=[n \pi]+1} \psi \pth{\frac{\varepsilon_t}{\sigma_0}}X_t=\frac{J_1(\psi)}{\sigma_0}\sum^{[n \pi^0]}_{t=[n \pi]+1} \varepsilon_t X_t+\sum^{[n \pi^0]}_{t=[n \pi]+1}  \sum_{q >1} \frac{J_q(\psi)}{q!} H_q(\nu_t) X_t:\equiv I_{1,n}+I_{2,n}
$$
where: $J_q(\psi)=\eE[\psi(\nu_1) H_q(\nu_1) ] $. On the other hand, as in the proof of  Theorem \ref{theorem 3}, we have  $I_{2,n}=o_{\eP}(I_{1,n})$. The  variance of $I_{1,n}$ is:
$$
\eE[I_{1,n}I_{1,n}^T]=\frac{J^2_1(\psi)}{\sigma^2_0} \sum^{[n(\pi^0-\pi)]}_{i=1}\sum^{[n(\pi^0-\pi)]}_{j=1} \gamma(|i-j|) \Gamma(|i-j|)$$
$$=[n(\pi^0-\pi)] \gamma(0) \Gamma(0)+2 \sum^{[n(\pi^0-\pi)]}_{i=1} [n(\pi^0-\pi)-i]
\gamma(i) \Gamma(i)$$
$$
=O\pth{L(n(\pi^0-\pi)){\cal L}^T(n(\pi^0-\pi))M(n(\pi^0-\pi)){\cal L} (n(\pi^0-\pi))}
$$
What implies:
$$
n^{-1}
\sum^{[n \pi^0]}_{t=[n \pi]+1} \psi \pth{\frac{\varepsilon_t}{\sigma_0}}X_t=O_{\eP} \pth{(n(\pi^0-\pi))^{-\min(\theta_i+\alpha)/2}L^{1/2}(n(\pi^0-\pi)){\cal L}^T(n(\pi^0-\pi)) {\cal L}(n(\pi^0-\pi))}
$$
This last relation together  with (\ref{e25}),  (\ref{e26}) and since the right-hand side of (\ref{e25}) is $O_{\eP}(n^{-k}\tilde L_1(n))$ imply: $
O_{\eP}(n^{-k}\tilde L_1(n))$
$$=O_{\eP}(n(\pi^0-\pi))+O_{\eP} \pth{(n(\pi^0-\pi))^{-\min(\theta_i+\alpha)/2}L^{1/2}(n(\pi^0-\pi)){\cal L}^T(n(\pi^0-\pi)) {\cal L}(n(\pi^0-\pi))}
$$
We obtain that: $\hat \pi_n-\pi^0=O_{\eP}(n^{-1-k} \tilde L_1(n))$.  \hspace*{\fill}$\blacksquare$ \\
\section{Lemmas}
\begin{lemma}
\label{lemma 1}
If solution $s_n(\xi) $ of equation (\ref{eq1}) exists, then it  is  well-defined,  bounded, strictly positive, with a probability arbitrarily large.
\end{lemma}
{\bf Proof of Lemma \ref{lemma 1}.} Since $\eE[r_t(\beta)=0]$ and $Var[r_t(\beta)]=Var[\varepsilon_t]+\beta Var[X] \beta^t <\infty$, by Bienaym\'e-Tchebichev inequality, we obtain that $r_t(\beta)$ is bounded with a probability arbitrarily large.\\ We prove that  $s_n(\xi)$ is bounded by reduction to absurdity. If $s_n(\xi)$ is not bounded then: there exists $\xi \in \Upsilon \times \Upsilon \times (0,1)$ and $n_\xi \in \N$ such that for all   $ n > n_\xi$, $ M>0$, exists  $ \epsilon >0$ such that: $\eP[s_n(\xi)>M] \geq 1-\epsilon$. Since $\rho$ is continuous and $\rho(0)=0$, then:
\begin{equation}
 \label{e2}
\rho\pth{\frac{r_t(\beta)}{s_n(\xi)}} \overset{\eP} {\underset{n \rightarrow \infty}{\longrightarrow}}0, \qquad t=1,...,n
\end{equation}
and
$$
\frac 1n \sum ^{[n\pi]}_{t=1}\rho \pth{\frac{r_t(\beta_1)}{s_n(\xi)}} +\frac 1n \sum ^{n}_{t=[n\pi]+1}\rho \pth{\frac{r_t(\beta_2)}{s_n(\xi)}}$$
 $$\leq \frac{[n\pi]}{n} \max_{1 \leq t \leq [n\pi]} \rho\pth{\frac{r_t(\beta_1)}{s_n(\xi)}}+\frac{n-[n\pi]}{n} \max_{[n\pi]+1 \leq t \leq n} \rho\pth{\frac{r_t(\beta_2)}{s_n(\xi)}}
$$
which, by (\ref{e2}), converges to 0 in probability, for $n\rightarrow \infty$. What is contradictory with (\ref{eq1}). To prove that $s_n(\xi)>0$, let us consider  function $g(\beta,s)=(\varepsilon-X\beta)/s$, with  $\beta$ in a  compact of $\R^d$ containing 0 and $s \in (0,\infty)$. Since $\varepsilon-X\beta$ is bounded with a    probability  close to  1, if $s_n(\xi)=0$, thus $\lim_{s \rightarrow 0} |g(\beta,s)|=\infty $, what is contradictory with  (\ref{eq1}). Hence, for all  $\epsilon>0$, there exists $\delta >0$ such that $\eP[\inf_{ \xi \in \Upsilon \times \Upsilon \times [0,1]} s_n(\xi) >\delta] >1-\epsilon$.
\hspace*{\fill}$\blacksquare$ \\

\begin{lemma}
 \label{lemma lDn} Under assumptions (A1)-(A3), for any  $\epsilon \in (0,1)$, $\xi \in \Upsilon \times \Upsilon \times [0,1]$, there exists a positive constant  $\delta$ such that: $\eP[{\underset{\xi \in \Upsilon \times \Upsilon \times [0,1]}{\inf}} D_n(\xi) >\delta] >1-\epsilon$.
\end{lemma}
{\bf Proof of Lemma \ref{lemma lDn}.}
Because $\xi$ belongs to a compact and taking into account relation  (\ref{H2}), we have to prove that for all $\epsilon>0$, $\xi \in \Upsilon \times \Upsilon \times [0,1]$, there exists  a $\delta>0$  such that:
$$
\eP\cro{n^{-1}\left[ \sum ^{[n\pi]}_{t=1}r_t(\beta_1) \psi\pth{\frac{r_t(\beta_1)}{s_n(\xi)}} +\sum ^{n}_{t=[n\pi]+1}r_t(\beta_2) \psi\pth{\frac{r_t(\beta_2)}{s_n(\xi)}} \right]>\delta}>1-\epsilon
$$
Since $r_t(\beta)$, $\psi\pth{\frac{r_t(\beta)}{s_n(\xi)}}$ have the same sign and since $\psi$ is continuous, we are going to show only that, for all  $\epsilon>0$, for all $\beta $ in compact set  $\Upsilon$, there exists a $\delta_1>0$ such that: $\eP[|\varepsilon-X \beta| >\delta_1]>1-\epsilon$.\\
Random variables $\varepsilon$ and $X$ are Gaussian and independent. Then: $\eP[|\varepsilon-X \beta| >\delta_1]=2 \eP[\varepsilon-X \beta <-\delta_1]=2 \Phi \pth{-\frac{\delta_1}{[\gamma(0)+\beta \Gamma(0) \beta^T]^{1/2}}}$. We recall that $\Phi$ denotes the standard Gaussian distribution. Then, the Lemma results by setting: $\delta_1=\inf_{\beta \in \Upsilon} [\gamma(0)+\beta \Gamma(0) \beta^T]^{1/2} \left| \Phi^{-1}\pth{\frac{1-\epsilon}{2}}\right|$ . 
\hspace*{\fill}$\blacksquare$ \\

The key for strong convergence proof is the following uniform convergence result.
\begin{lemma}
\label{lemma 2} For all $\varrho >0$, under assumptions (A1)-(A3),  for \\
$\Omega_\varrho(\xi)=\acc{\xi^* \in \Upsilon \times \Upsilon \times [0,1]; \| \eta-\eta^*\| <\varrho, |\pi-\pi^*| <\varrho}$, we have:
$$
\eE \cro{\sup_{\xi^* \in \Omega_\varrho(\xi)} |s_n(\eta,\pi) - s_n(\eta^*,\pi^*)|} {\underset{\varrho \rightarrow 0}{\longrightarrow}}0
$$ 

\end{lemma}
{\bf Proof of Lemma \ref{lemma 2}.}
We have the triangular inequality:\\
$|s_n(\eta,\pi) -s_n(\eta^*,\pi^*)| \leq |s_n(\eta,\pi) -s_n(\eta,\pi^*)|+|s_n(\eta,\pi^*)-s_n(\eta^*,\pi^*)|$. First, we will study $s_n(\eta,\pi^*) -s_n(\eta^*,\pi^*)$. By the mean value theorem  (TVM), we have:
\begin{equation}
 \label{e7}
s_n(\eta,\pi^*)-s_n(\eta^*,\pi^*)=(\beta_1-\beta^*_1) \frac{\partial s_n}{\partial \beta_1}(\tilde \beta_1,\beta^*_2,\pi^* )+(\beta_2-\beta^*_2)\frac{\partial s_n}{\partial \beta_2}(\beta^*_1,\tilde \beta_2,\pi^* )
\end{equation}
where $\tilde \beta_1=\beta_1+\upsilon_1(\beta_1-\beta^*_1)$, $\tilde \beta_2=\beta_2+\upsilon_2(\beta_2-\beta^*_2)$, $\upsilon_1,\upsilon_2 \in (0,1)$.
 By Lemma \ref{lemma lDn}, applying  Cauchy-Schwarz inequality in (\ref{e6}) and taking into account that $\psi$ is bounded, we obtain:
\begin{equation}
\label{eSC}
\eE \cro{\left| \frac{\partial s_n(\tilde \beta_1,\beta^*_2,\pi^* )}{\partial \beta_1} \right|} \leq Cn^{-1} \sum ^{[n\pi]}_{t=1} \pth{\eE [X_t^2] \eE\left[ \psi^2 \pth{\frac{r_t(\tilde \beta_1)}{s_n(\tilde \beta_1,\beta^*_2,\pi^* )}}\right]}^{1/2}  <C
\end{equation}
Then, writing a similar relation for $({\partial s_n}/{\partial \beta_2})(\beta^*_1,\tilde \beta_2,\pi^* )$, we have for (\ref{e7}):
\begin{equation}
 \label{e8}
\eE \cro{|s_n(\eta,\pi^*)-s_n(\eta^*,\pi^*)|} \longrightarrow 0, \qquad \textrm{for } \varrho \rightarrow 0
\end{equation}
Let us remark that if $\pi^*=0$ or $\pi^*=1$, then in relation (\ref{e7}), the term in $\beta_1$, respectively $\beta_2$, does not appear.\\ 
Now, we study $|s_n(\eta,\pi) -s_n(\eta,\pi^*)|$, supposing that $\pi < \pi^*$. Since $s_n(\eta,\pi)$ and  $s_n(\eta,\pi^*)$  are both solutions of (\ref{eq1}), we have:
$$
n^{-1} \sum ^{[n\pi^*]}_{t=1}\cro{\rho \pth{\frac{r_t(\beta_1)}{s_n(\eta,\pi)}} - \rho \pth{\frac{r_t(\beta_1)}{s_n(\eta,\pi^*)}}}+ n^{-1} \sum^{n}_{t=[n\pi^*]+1} \cro{ \rho \pth{\frac{r_t(\beta_2)}{s_n(\eta,\pi)}} - \rho \pth{\frac{r_t(\beta_2)}{s_n(\eta,\pi^*)}} }$$
$$
=n^{-1} \sum^{[n\pi^*]}_{[n\pi]+1} \cro{\rho \pth{\frac{r_t(\beta_1)}{s_n(\eta,\pi)}} -\rho \pth{\frac{r_t(\beta_2)}{s_n(\eta,\pi)}}}
$$
Thus, applying the MVT:
 \begin{equation}
\label{e10}
\begin{array}{ll}
 n^{-1} \cro{s_n(\eta,\pi^*)-s_n(\eta,\pi)} \cro{\sum ^{[n\pi^*]}_{t=1} r_t(\beta_1)\psi \pth{\frac{r_t(\beta_1)}{u^{(1)}_n(\eta,\pi,\pi^*)}}+\sum^{n}_{t=[n\pi^*]+1}r_t(\beta_2)\psi \pth{\frac{r_t(\beta_2)}{u^{(2)}_n(\eta,\pi,\pi^*)}} }\\
=n^{-1} \sum^{[n\pi^*]}_{[n\pi]+1}  [X_t(\beta_1-\beta_2)] \psi \pth{\frac{\tilde r_t(\beta_1,\beta_2)}{s_n(\eta,\pi)}}
\end{array}
\end{equation}
where $u^{(1)}, u^{(2)}$ are two positive bounded  functions, not necessarily solutions of  (\ref{eq1}) and $\tilde r_t(\beta_1,\beta_2)=r_t(\beta_1)+m_t [r_t(\beta_2)-r_t(\beta_1)]$, with $0 < m_t <1$.
By relation (\ref{H2}):
\begin{equation}
 \label{e9}
\sum ^{[n\pi^*]}_{t=1} r_t(\beta_1)\psi \pth{\frac{r_t(\beta_1)}{u^{(1)}_n(\eta,\pi,\pi^*)}}+\sum^{n}_{t=[n\pi^*]+1}r_t(\beta_2)\psi \pth{\frac{r_t(\beta_2)}{u^{(2)}_n(\eta,\pi,\pi^*)}} >0
\end{equation}
with a probability close to 1. On the other hand: $\tilde r_t(\beta_1,\beta_2)=Y_t-X_t[\beta_1+m_t(\beta_2-\beta_1)]=r_t(\beta_1+m_t(\beta_2-\beta_1))$.
Using the same arguments as for (\ref{eSC}),  we obtain that:
$$
n^{-1} \sum ^{[n\pi^*]}_{t=[n \pi]+1}\pth{\eE[ X_t^2] \eE\left[ \psi^2 \pth{\frac{\tilde r_t(\beta_1,\beta_2)}{s_n(\eta,\pi)}}\right] }^{1/2} \leq C_1 (\pi^*-\pi)
$$
where $C_1$ is a  vector with all bounded components. Taking into account also (\ref{e9}), we obtain for (\ref{e10}):
\begin{equation}
 \label{e11}
\eE[|s_n(\eta,\pi) -s_n(\eta,\pi^*)|] \leq C \| \beta_1- \beta_2\| \cdot |\pi-\pi^*| <C \varrho {\underset{\varrho \rightarrow 0}{\longrightarrow}}0
\end{equation}
Relations (\ref{e8}) and (\ref{e11}) imply the Lemma.
\hspace*{\fill}$\blacksquare$ \\


\begin{thebibliography}{00}
\bibitem{Bai:94} J. Bai, 
\newblock Least squares estimation of a shift in linear processes,
\newblock { Journal of Time Series Analysis}  15 (1994) 453-472.
\bibitem{Bai:98}  J. Bai,  
\newblock Estimation of multiple-regime regressions with least absolute deviation,
\newblock { Journal of Statistical Planning Inference},  74 (1998)  103-134.
\bibitem{Bai:Perron:98} J. Bai, P.  Perron, 
\newblock Estimating and testing linear models with multiple structural changes,
\newblock { Econometrica} 66  (1998) 47-78.
\bibitem{Baillie:96}  R. Baillie,  
\newblock Long memory processes and fractional integration in econometrics,
\newblock { Journal of Econometrics},  73 (1996)  5-59.
\bibitem{Beaton:Tukey:74} A.E. Beaton, J.W. Tukey, 
\newblock The fitting of power series, meaning polynomials, illustrated on band-spectroscopic data,
\newblock { Technometrics} 16 (1974)  147-185.
\bibitem{Beran:94} J. Beran,  
\newblock  Statistics for long-memory process, Chapman $\&$ Hall, New York, 1994.
\bibitem{Bhattacharya:94} P.K.Bhattacharya,    
\newblock  { Some aspects of change-point analysis}.
\newblock {IMS Lecture Notes-Monograph Series}, Vol. 23,   Hayward, CA 1994, pp.28-56.
\bibitem{Cheung:93} Y.X. Cheung,  
\newblock  Long memory in foreign exchange rates, Journal of Bussiness and Economic Statistics 11(1993) 93-101.
\bibitem{Ciuperca:Dapzol:08} G. Ciuperca, N. Dapzol, 
\newblock  Maximum likelihood estimator in a multi-phase random regression model,
\newblock { Statistics}  42 (2008) 363-381.
\bibitem{Ciuperca:09} G. Ciuperca 
\newblock Estimating nonlinear regression with and without change-points by the LAD-method,
\newblock { Annals of the Institute of Statistical Mathematics} in revision.
\bibitem{Davies:90} L. Davies, 
\newblock The Asymptotics of S-Estimators in the Linear Regression Model,
\newblock { Annals of Statistics} 18 (1990) 1651-1675.
\bibitem{Ding:Granger:Engle:93} Z. Ding,  C.W.J. Granger, R.F.  Engle, 
\newblock  A long memory property of stock market returns and a new model,
\newblock { Journal of Empirical Finance } 1 (1993) 83-106.
\bibitem{Feder:75a} P.I. Feder,  
\newblock  On asymptotic distribution theory in segmented regression problems-identified case,
\newblock { Ann. Statist.}  3 (1975) 49-83.
\bibitem{Feder:75b} P.I. Feder,  
\newblock  The log likelihood ratio  in segmented regression,
\newblock { Ann. Statist.}  3 (1975) 84-97.
\bibitem{Fiteni:02} I. Fiteni, 
\newblock Robust estimation of structural break points,
\newblock { Econometric Theory} 18 (2002) 349-386.
\bibitem{Fiteni:04} I. Fiteni, 
\newblock $\tau$-estimators of regression models with structural change of unknown location,
\newblock { Journal of Econometrics} 119 (2004) 19-44.
\bibitem{Guo: Koul:07} H. Guo, H.L. Koul,  
\newblock Nonparametric regression with heteroscedastic long memory errors,
\newblock { Journal of Statistical Planning Inference} 137 (2007) 379-404.
\bibitem{Hampel:71} F.R.. Hampel, 
\newblock A general quantitative definition of robustness,
\newblock { Annals of  Mathematical  Statistics} 42 (1971) 1887-1896.
\bibitem{Hidalgo:Robinson:96} J. Hidalgo, P.M. Robinson,  
\newblock Testing for structural change in a long-memory environment,
\newblock { Journal of Econometrics} 1 (1996) 159-174.
\bibitem{Horvath: Kokoszka:97} L. Horvath, P. Kokoszka,  
\newblock The effect of long-range dependence on change-point estimators,
\newblock { Journal of Statistical Planning Inference} 64 (1997) 57-81.
\bibitem{Huber:67} P.J. Huber,   
The behaviour of maximum likelihood estimates under nonstandard conditions. Proceedings of the Fifth Berkeley Symposium on Mathematics Statistic and Probability, Vol 1, University California Press, Berkeley, 1967, pp  221-234.  
\bibitem{Kim:Kim:08} J.  Kim , H.J Kim,   
\newblock Asymptotic results in segmented multiple regression,
\newblock { Journal of Multivariate Analysis} 99 (2008) 2016-2038.
\bibitem{Koul:Qian:02} H.L. Koul, L.  Qian,  
\newblock  Asymptotics of maximum likelihood estimator in a two-phase linear regression model,
\newblock { Journal of Statistical Planning and Inference}   108 (2002) 99-119.
\bibitem{Koul:Baillie:03} H.L. Koul, R.T. Baillie,  
\newblock Asymptotics of M-estimators in non-linear regression with long memory design,
\newblock {Statistics and Probability Letters} 61 (2003) 237-252.
\bibitem{KQS} H.L. Koul,  L. Qian, D.  Surgailis, 
\newblock  Asymptotics of M-estimators in two-phase linear regression models,
\newblock { Stochastic Processes and their Applications } 103 (2003) 123-154.
\bibitem{Kuan:Hsu:98} C.M. Kuan, C.C.  Hsu,  
\newblock  Change-point estimation of fractionally integrated processes,
\newblock { Journal of Time Series Analysis}   19 (1998) 693-708.
\bibitem{Lazarova:05} S. Lazarov\'a, 
\newblock Testing for structural change in regression with long memory processes,
\newblock {Journal of  Econometrics} 129 (2005) 329-372.
\bibitem{Lo:91} A.W. Lo, 
\newblock Long term memory in stock market prices,
\newblock { Econometrica} 59 (1991) 1279-1313.
\bibitem{Nunes:Kuan:Newbold:95} L.C. Nunes, C.M. Kuan,  P. Newbold, 
\newblock  Spurious break
\newblock {Econometric Theory} 11(1995) 736-749.
\bibitem{Palma:07}  W. Palma, 
\newblock Long-Memory Time Series, Theory and Methods, 
\newblock Wiley, New Jersey, 2007.
\bibitem{Robinson:94}  P.M. Robinson, 
\newblock Time series with strong dependence,
\newblock In: C.A. Sims, Advances in Econometrics, Sixth World Congress, Cambridge Univ. Press, 1994.
\bibitem{Roelant:VanAelst:Croux:08} E. Roelant, S. Van Aelst,  C. Croux, 
\newblock  Multivariate generalized S-estimators
\newblock { Journal of Multivariate Analysis} (2008) in press.
\bibitem{Rousseeuw:Yohai:84} P.J. Rousseeuw,    V.J. Yohai,   
\newblock Robust regression by means of S-estimators,
\newblock In: Robust and Nonlinear Time Series Analysis, Franke, J., Härdle, W., Martin, R.D. (Eds.), Lecture Notes in Statistics, Vol. 26, Springer, New York, 1984, pp. 256-–272.
\bibitem{Rukhin:Vajda:97} A.L. Rukhin,  I. Vajda,  
 \newblock Change-point estimation as a nonlinear regression problem,
\newblock { Statistics} 30 (1997) 181-200.
\bibitem{Sibbertsen:04}  Ph. Sibbertsen, 
\newblock Long memory versus structural breaks: an overview,
\newblock { Statistical Paper} 4 (2004) 465-515.
\bibitem{Zhengyan:Degui:Jia:07} L. Zhengyan,L. Degui, C. Jia, 
\newblock  Asymptotic behavior for $S$-estimators in random design linear model with long-range-dependent errors,
\newblock { Metrika} 66 (2007) 289-303.

\end{thebibliography}
\end{document}